\documentclass[review,11pt]{article}
\usepackage{amsmath}
\usepackage{amssymb}
\usepackage{amsthm}
\usepackage{enumerate}
\usepackage[mathscr]{eucal}
\usepackage{eqlist}
\usepackage[dvips]{graphicx}
\usepackage{subfigure}
\usepackage{stmaryrd}
\usepackage{amsfonts}
\usepackage{}
\usepackage{mathrsfs}
\usepackage{amssymb}
\usepackage{latexsym,amssymb,amsmath,amsfonts,amsthm}

\makeatletter
\newif\ifmsbmloaded@
\def\loadmsbm{\msbmloaded@true
  \font\tenmsb=msbm10 scaled 1\@ptsize00
  \font\sevenmsb=msbm7 scaled 1\@ptsize00
  \font\fivemsb=msbm5 scaled 1\@ptsize00
  \alloc@8\fam\chardef\sixt@@n\msbfam
  \textfont\msbfam=\tenmsb
  \scriptfont\msbfam=\sevenmsb
  \scriptscriptfont\msbfam=\fivemsb
  }
\def\nonmatherr@#1{\errmessage%
{LateX error: \string#1\space allowed only in math mode}}
\def\Bbb{\relax\ifmmode\expandafter\Bbb@\else
  \expandafter\nonmatherr@\expandafter\Bbb\fi}
\def\Bbb@#1{{\Bbb@@{#1}}}
\def\Bbb@@#1{\fam\msbfam\relax#1}
\@addtoreset{equation}{section}

\makeatother \loadmsbm

\newcommand{\beq}{\begin{equation}}
\newcommand{\eeq}{\end{equation}}
\newcommand{\ben}{\begin{eqnarray}}
\newcommand{\een}{\end{eqnarray}}
\newcommand{\beno}{\begin{eqnarray*}}
\newcommand{\eeno}{\end{eqnarray*}}

\begin{document}
\title{Grand Fujii-Fujii-Nakamoto operator inequality dealing with operator order and operator chaotic order}

\author{ Jian Shi$^{a}$*\\
     {\small  $^{a}$College of Mathematics and Information Science, Hebei University,}\\
     {\small  Baoding, 071002, P.R. China.}\\
     {\small * Corresponding author, E-mail: mathematic@126.com }\\
     }
\date{}
\maketitle

\vspace{-1.2in} \vspace{.9in} \vspace{0.2cm} {\bf Abstract.} In this paper, we shall prove that a grand Fujii-Fujii-Nakamoto operator inequality implies operator order and operator chaotic order under different conditions.

\vspace{0.2cm} {\bf Key words.}  operator order, operator chaotic order, Furuta inequality.
 \vspace{0.2cm}


\vskip0.2cm

\section{Introduction}

\ \ \ \ A capital letter, such as $T$, stands for a bounded linear operator on a Hilbert space $\mathcal{H}$.
$T\geq 0$ means $T$ is a positive operator and $T>0$ means $T$ is a positive and invertible operator.

$A\geq B$ and $\log A\geq \log B$ are called operator order and operator chaotic order, respectively.

There are many sufficient conditions and necessary conditions of operator order and operator chaotic order, such as\cite{Ando1987}, \cite{Fujii2014}, \cite{Fujii1993}, \cite{Yamazaki2002}.

In 1993, M. Fujii, T. Furuta and E. Kamei obtained the following result on operator chaotic order.

\noindent{\bf Theorem 1.1}(\cite{Fujii1993}). If $A, B>0$ satisfy
$$
(B^{\frac 1 2}A^{\frac {1-t}{2}}B^{t}A^{\frac {1-t}{2}}B^{\frac 1 2})^{\frac 1 2}\geq B \eqno (1.1)
$$
for all $t>1$, then $\log B\geq \log A$.

In 2014, J. I. Fujii, M. Fujii and R. Nakamoto proposed an extension of Theorem 1.1, by which they make use of a beautiful operator inequality(we call it Fujii-Fujii-Nakamoto operator inequality here).

\noindent{\bf Theorem 1.2}(\cite{Fujii2014}). If $A, B>0$ satisfy
$$
(B^{\frac s 2}A^{\frac {s-t}{2}}B^{t}A^{\frac {s-t}{2}}B^{\frac s 2})^{\frac {1}{2s}}\geq B \eqno (1.2)
$$
for some $t>s>0$, then the following result hold:\\
(I) If $t\geq 3s-2\geq 0$, then $\log B\geq \log A$, and if $t\geq s+2$ is additionally assumed, then $B\geq A$;\\
(II) If $0<s<{\frac 1 2}$, then $\log B\geq \log A$, and if $t\geq s+2$ is additionally assumed, then $B\geq A$.

In this paper, we will extend  Theorem 1.2 and show some results on operator order and operator chaotic order.

In order to prove our results, we list a famous operator inequality here, which is called Furuta inequality.

\noindent{\bf Theorem 1.3} (\cite{Furuta1987}.) If $A\geq B\geq 0$, then
$$
(A^{\frac r 2}A^{p}A^{\frac r 2})^{\frac 1 q}\geq (A^{\frac r 2}B^{p}A^{\frac r 2})^{\frac 1 q}
$$
and
$$
(B^{\frac r 2}A^{p}B^{\frac r 2})^{\frac 1 q}\geq (B^{\frac r 2}B^{p}B^{\frac r 2})^{\frac 1 q}
$$
hold for $r\geq 0, p\geq 0, q\geq 1$ with $(1+r)q\geq p+r.$

\section{Main Results}

\ \ \ \ In this section, we shall show the main results and prove them.

\noindent{\bf Theorem 2.1.} If $A, B>0, $  satisfying
$$( \underbrace{B^{\frac s n}A^{\frac {s-t}{n}}B^{\frac {2t}{n}}A^{\frac {s-t}{n}}B^{\frac {2t}{n}}\cdots B^{\frac {2t}{n}}A^{\frac {s-t}{n}}  B^{\frac {2t}{n}}A^{\frac {s-t}{n}}B^{\frac s n}}_{ \ A^{\frac {s-t}{n}}\ appears\ n\ times, \ B^{\frac {2t}{n}}  appears\ n-1 \ times  } )^{\frac {n}{(n+2)s+(n-2)t}}\geq B \eqno (2.1)$$
for some $t>s>0$ and $n\geq 2$, then the following results hold:\\
(I) If $(n+2)s+(n-2)t>n\geq 3s-t$, then $\log B\geq \log A$, and if additional condition $t-s\geq n$ is assumed, then $B\geq A$;\\
(II) If $(n+2)s+(n-2)t\leq n$, then $\log B\geq \log A$, and if additional condition $t-s\geq n$ is assumed, then $B\geq A$.

\noindent{\bf Proof of (I).} Firstly, put $A_{1}=( \underbrace{B^{\frac s n}A^{\frac {s-t}{n}}B^{\frac {2t}{n}}A^{\frac {s-t}{n}}B^{\frac {2t}{n}}\cdots B^{\frac {2t}{n}}A^{\frac {s-t}{n}}  B^{\frac {2t}{n}}A^{\frac {s-t}{n}}B^{\frac s n}}_{ \ A^{\frac {s-t}{n}}\ appears\ n\ times, \ B^{\frac {2t}{n}}  appears\ n-1 \ times  } )^{\frac {n}{(n+2)s+(n-2)t}}$, $l={\frac{(n+2)s+(n-2)t} {n}}$. Obviously, $l>1$, $2{\frac {t-s}{n}}>0$ and conditions in Furuta inequality are satisfied.

It follows that
$$
(B^{\frac {t-s}{n}}A_{1}^{l}B^{\frac {t-s}{n}})^{\frac {1+2{\frac {t-s}{n}}}{l+2{\frac {t-s}{n}}}}\geq B^{1+2{\frac {t-s}{n}}}. \eqno (2.1)
$$
Then we have
\begin{equation}
\begin{split}
&\ ( \underbrace{B^{\frac t n}A^{\frac {s-t}{n}}B^{\frac {2t}{n}}A^{\frac {s-t}{n}}B^{\frac {2t}{n}}\cdots B^{\frac {2t}{n}}A^{\frac {s-t}{n}}  B^{\frac {2t}{n}}A^{\frac {s-t}{n}}B^{\frac t n}}_{ \ A^{\frac {s-t}{n}}\ appears\ n\ times, \ B^{\frac {2t}{n}}  appears\ n-1 \ times  } )^{\frac {n+2(t-s)}{nl+2(t-s)}}\geq B^{\frac {n+2(t-s)}{n}}\\
=&  (B^{\frac t n}A^{\frac {s-t}{n}}B^{\frac { t}{n}})^{n{\frac {n+2(t-s)}{nl+2(t-s)}}}
\end{split}\tag{2.2}
\end{equation}

Notice that ${n{\frac {n+2(t-s)}{nl+2(t-s)}}} \geq 1$ due to the fact that $n\geq 3s-t$. Apply L\"{o}wner-Heinz inequality to (2.2), then the following inequality is obtained.
$$
B^{\frac t n}A^{\frac {s-t}{n}}B^{\frac { t}{n}}\geq B^{\frac {s+t}{n}}, \eqno (2.3)
$$
which is equivalent to
$$
A^{\frac {s-t}{n}}\geq B^{\frac {s-t}{n}}. \eqno (2.4)
$$

Taking reverse and then taking logarithm on both sides of (2.4), $\log B\geq \log A$ is obtained.

If $t-s\geq n$, taking reverse and applying L\"{o}wner-Heinz inequality to $0<{\frac {n}{t-s}}\leq 1$, $B\geq A$ holds, obviously.

\noindent{\bf Proof of (II).}
According to L\"{o}wner-Heinz inequality, (2.1) implies that
$$
 \underbrace{B^{\frac s n}A^{\frac {s-t}{n}}B^{\frac {2t}{n}}A^{\frac {s-t}{n}}B^{\frac {2t}{n}}\cdots B^{\frac {2t}{n}}A^{\frac {s-t}{n}}  B^{\frac {2t}{n}}A^{\frac {s-t}{n}}B^{\frac s n}}_{ \ A^{\frac {s-t}{n}}\ appears\ n\ times, \ B^{\frac {2t}{n}}  appears\ n-1 \ times  }
\geq B^{\frac {(n+2)s+(n-2)t}{n}}.  \eqno (2.5)
$$
It follows that
\begin{equation}
\begin{split}
&\ \underbrace{B^{\frac t n}A^{\frac {s-t}{n}}B^{\frac {2t}{n}}A^{\frac {s-t}{n}}B^{\frac {2t}{n}}\cdots B^{\frac {2t}{n}}A^{\frac {s-t}{n}}  B^{\frac {2t}{n}}A^{\frac {s-t}{n}}B^{\frac t n}}_{ \ A^{\frac {s-t}{n}}\ appears\ n\ times, \ B^{\frac {2t}{n}}  appears\ n-1 \ times  }
\geq B^{s+t}\\
=&  (B^{\frac t n}A^{\frac {s-t}{n}}B^{\frac { t}{n}})^{n}
\end{split}\tag{2.6}
\end{equation}
Then we have
$$B^{\frac t n}A^{\frac {s-t}{n}}B^{\frac { t}{n}}\geq B^{\frac {s+t}{n}}\eqno (2.7)$$ by L\"{o}wner-Heinz inequality.

(2.7) is equivalent to $A^{\frac {s-t}{n}}\geq B^{\frac {s-t}{n}}$.

By the same discussion as in proof in (I), we complete the proof.\quad \quad   $\Box$\\

\noindent{\bf Corollary 2.1.} If $A, B>0, $  satisfying
$$(  {B^{\frac s 3}A^{\frac {s-t}{3}}B^{\frac {2t}{3}}A^{\frac {s-t}{3}}B^{\frac {2t}{3}} A^{\frac {s-t}{3}}   B^{\frac s 3}}  )^{\frac {3}{ 5s+t}}\geq B \eqno (2.8) $$
for some $t>s>0$, then the following results hold:\\
(I) If $5s+t>3\geq 3s-t$, then $\log B\geq \log A$, and if additional condition $t-s\geq 3$ is assumed, then $B\geq A$;\\
(II) If $5s+t\leq 3$, then $\log B\geq \log A$.

\noindent{\bf Proof. }  We only need to set $n=3$ in Theorem 2.1.\\

\noindent{\bf Remark 2.1.} If $n=2$, Theorem 2.1 is just Theorem 1.2.\\

Similarly, we can obtain another generalized Fujii-Fujii-Nakamoto operator inequality which extend Theorem 1.2.\\
\noindent{\bf Theorem 2.2}. If $A, B>0$ satisfy
$$
(B^{\frac s p}A^{\frac {s-t}{p}}B^{\frac {2t}{p}}A^{\frac {s-t}{p}}B^{\frac s p})^{\frac {p}{4s}}\geq B \eqno (2.9)
$$
for some $t>s>0$, then the following result hold:\\
(I) If $4s>p\geq 3s-t$, then $\log B\geq \log A$, and if $t-s\geq p$ is additionally assumed, then $B\geq A$;\\
(II) If $p\geq 4s$, then $\log B\geq \log A$, and if $t-s\geq p$ is additionally assumed, then $B\geq A$.

\noindent{\bf Remark 2.2.} The proof of Theorem 2.2 is the same as in \cite{Fujii2014}, and we omit it here. If $p=2$, Theorem 2.2 is just Theorem 1.2.\\

\noindent{\bf Corollary 2.2.}  If $A, B>0$ satisfy
$$
(B^{\frac s 3}A^{\frac {s-t}{3}}B^{\frac {2t}{3}}A^{\frac {s-t}{3}}B^{\frac s 3})^{\frac {3}{4s}}\geq B \eqno (2.10)
$$
for some $t>s>0$, then the following result hold:\\
(I) If $4s>3\geq 3s-t$, then $\log B\geq \log A$, and if $t-s\geq 3$ is additionally assumed, then $B\geq A$;\\
(II) If $3\geq 4s$, then $\log B\geq \log A$, and if $t-s\geq 3$ is additionally assumed, then $B\geq A$.\\

\noindent{\bf Proof. }  We only need to set $p=3$ in Theorem 2.2.\\

\noindent {\bf Acknowledgements.}   J. Shi (the corresponding author) is supported by National Natural Science Foundation of China(No. 11702078 and No. 61702019), Hebei Education Department (No. ZC2016009).

\end{document}